# Generalized Estimators Using Characteristics of Poisson distribution

Prayas Sharma, Hemant K. Verma, Nitesh K. Adichwal and *Rajesh Singh

Department of Statistics, Banaras Hindu University

Varanasi(U.P.), India-221005

* Corresponding author

rsinghstat@gmail.com

**Abstract**

In this article, we have proposed a generalized class of estimators, exponential class of estimators based on adaption of Sharma and Singh (2015) and Solanki and Singh (2013) and simple difference estimator for estimating unknown population mean in case of Poisson distributed population in simple random sampling without replacement. The expressions for mean square errors of the proposed classes of estimators are derived to the first order of approximation. It is shown that the adapted version of Solanki and Singh (2013), exponential class of estimator, is always more efficient than usual estimator, ratio, product, exponential ratio and exponential product type estimators and equal efficient to simple difference estimator. Moreover, the adapted version of Sharma and Singh (2015) estimator are always more efficient than all the estimators available in literature. In addition, theoretical findings are supported by an empirical study to show the superiority of the constructed estimators over others with earthquake data of turkey.

## 1. Introduction

In the sampling literature, it is well known that efficiency of the estimator of population parameters of a study variable y can be increased by the use of auxiliary information related to x, which is highly correlated with study variable y. Several authors including Singh and Kumar (2011), Singh and Solanki (2012) and Sharma and Singh(2014a,b) suggested estimators using auxiliary information under different situations but without considering the distributions of study and auxiliary variate. However some distribution, like Poisson distribution is generally used for the natural population to express the probability of a number of rare events. Aftershocks constitute the greatest proportions of shocks in an earthquake catalogue and if aftershocks are effectively considered, they can give us information for understanding the whole cycle of seismic activity. Due to this reason estimating the number of aftershocks is important in seismology and has received much attention in resent literature (see Ozel,2011). Ratio type estimators can give us information about the number of aftershocks in a specified region. Consequently, in this paper we have proposed a general class of exponential estimators and difference estimator for population mean using auxiliary information from a Poisson distributed population. Here we used the earthquake data for empirical study since earthquake are rare events and generally follows a Poisson distribution.

Consider a population $U = (u_1, u_2, ..., u_n)$ of size N identifiable and distinct units. Let y and x be the study and auxiliary variables associated with each unit $u_j = (j = 1, 2, ...... N)$ of the population respectively. Assume that X's are known units and Y's are unknown units for all the population. Suppose a random sample of size n is drawn using simple random sampling without replacement (SRSWOR) from the population. Let us assume that the parent population has a Poisson distribution. We know that the nature of the sampling distribution depends on the nature of the population from which the random sample is drawn therefore, random samples which are drawn from a Poisson distrusted population follows also a Poisson distribution. Further les us select observations $(y_i, x_i)$, (i=1,2.......,n) from a Poisson distributed population. Using this sampling design we can define classical ratio estimator as

$$t_{ratio} = \overline{y}_{po} \left( \frac{\overline{X}}{\overline{x}_{po}} \right) \qquad (1.1)$$

where $\overline{y}_{po}$ and $\overline{x}_{po}$ are the sample means of the study and auxiliary variables from Poisson distributed population, respectively . if we suppose that auxiliary variable x has Poisson distribution with parameter $\lambda_1 > 0$, then the values of the parameters of auxiliary variable x are given by $\overline{X} = \lambda_1$, $S_x = \sqrt{\lambda_1}$, $C_x = \frac{S_x}{\overline{X}} = \frac{1}{\sqrt{\lambda_1}}$ respectively, further let the study variable y has a Poisson distribution with parameter $\lambda_2 > 0$, then the expression for study variable can be written as $\overline{Y} = \lambda_2$, $S_y = \sqrt{\lambda_2}$, $C_y = \frac{S_y}{\overline{Y}} = \frac{1}{\sqrt{\lambda_2}}$. The correlation coefficient between the study variable y and auxiliary variable x is obtained by using trivariate reduction method. The trivariate reduction method is an appealing method for constructing bivariate Poisson distribution (See Lai (1995)). The method is to create a pair of dependent Poisson distributed random variable from three independent Poisson distributed random variables. Let k, w and z are independently Poisson distributed observations from a Poisson distributed population then a bivariate Poisson distribution of study variable y and auxiliary variable x is generated by setting $x_i = k_i + z_i$ and $y_i = w_i + z_i$ for i=1,2,.....,n. Assuming that the parameters of k, w and z are $\gamma_1, \gamma_2$ and $\gamma_3$ respectively. The correlation coefficient between the study variable y and auxiliary variable x is defined as

$$\rho_{x_{po}y_{po}} = \frac{Cov(x,y)}{\sqrt{S_x^2 S_y^2}} = \frac{[(\gamma_1+\gamma_3)(\gamma_2+\gamma_3)+\gamma_3]-(\gamma_1+\gamma_3)(\gamma_2+\gamma_3)}{\sqrt{(\gamma_1+\gamma_3)(\gamma_2+\gamma_3)}} = \frac{\gamma_3}{\sqrt{(\gamma_1+\gamma_3)(\gamma_2+\gamma_3)}}$$

Here, $\rho_{x_{po}y_{po}}$ is restricted to be strictly positive since $\gamma_1, \gamma_2$ and $\gamma_3$ are always positive. Since we have drawn the observations $(y_i, x_i)$, i=1,2,...,n from a Poisson distributed population with parameters $\lambda_1$ and $\lambda_2$, then we have $\bar{x}_{po} = \sum_{i=1}^{n}(k_i + z_i)/n$ and $\bar{y}_{po} = \sum_{i=1}^{n}(w_i + z_i)/n$. The covariance between $\bar{x}_{po}$ and $\bar{y}_{po}$ is defined as

$$Cov(\bar{x}_{po}, \bar{y}_{po}) = E(\bar{x}_{po}, \bar{y}_{po}) - E(\bar{x}_{po})E(\bar{y}_{po})$$

$$\left[\frac{(\gamma_1+\gamma_3)}{n}\frac{(\gamma_2+\gamma_3)}{n} + \frac{\gamma_3}{n}\right] - \frac{(\gamma_1+\gamma_3)}{n}\frac{(\gamma_2+\gamma_3)}{n} = \frac{\gamma_3}{n}$$

To fine the bias and mean square error, let us define

$$e_0 = \frac{\bar{y}_{po} - \bar{Y}}{\bar{Y}}, \quad e_1 = \frac{\bar{x}_{po} - \bar{X}}{\bar{X}}$$

$$E(e_0^2) = \frac{V(\bar{y}_{po})}{\bar{Y}^2} = \frac{\lambda_2}{n} = \frac{(\gamma_2+\gamma_3)}{n}, \quad E(e_1^2) = \frac{V(\bar{x}_{po})}{\bar{X}^2} = \frac{\lambda_1}{n} = \frac{(\gamma_1+\gamma_3)}{n}$$

$$E(e_0 e_1) = \frac{Cov(\bar{x}_{po}, \bar{y}_{po})}{\bar{X}\bar{Y}} = \frac{\gamma_2}{n}$$

Now expressing estimator $t_{ratio}$ in term of e's, we have

$$(t_r - \bar{Y}) = \bar{Y}(-e_1 + e_1^2 + e_0 - e_0 e_1)$$

The Bias and MSE of the estimator $t_{ratio}$ are respectively, given by

$$\text{Bias}(t_{ratio}) = \overline{Y}\left(\frac{\gamma_1}{n}\right) \qquad (1.2)$$

$$\text{MSE}(t_{ratio}) = \overline{Y}^2\left(\frac{\gamma_1 + \gamma_2}{n}\right) \qquad (1.3)$$

## 2. Estimators in Literature

Koyuncu and Ozel (2013) suggested an exponential ratio estimator for estimation unknown population mean as

$$t_{k1} = \overline{y}_{po}\left(\frac{\overline{X} - \overline{x}_{po}}{\overline{X} + \overline{x}_{po}}\right) \qquad (2.1)$$

The bias and MSE of the estimator $t_{k1}$ are respectively, given by

$$\text{Bias}(t_{k1}) = \overline{Y}\left(\frac{3\gamma_1 - \gamma_3}{8n}\right) \qquad (2.2)$$

$$\text{MSE}(t_{k1}) = \overline{Y}^2\left(\frac{\gamma_1 + 4\gamma_2 + \gamma_3}{4n}\right) \qquad (2.3)$$

In case of negative correlation coefficient between study variable y and auxiliary variable x Koyuncu and Ozel (2013) suggested exponential product estimator as

$$t_{k2} = \overline{y}_{po}\left(\frac{\overline{x}_{po} - \overline{X}}{\overline{x}_{po} + \overline{X}}\right) \qquad (2.4)$$

The bias and MSE of the estimator $t_{k1}$ are respectively, given by

$$\text{Bias}(t_{k2}) = \overline{Y}\left(\frac{\gamma_1 + 5\gamma_3}{8n}\right) \qquad (2.5)$$

$$\text{MSE}(t_{k2}) = \overline{Y}^2\left(\frac{\gamma_1 + 4\gamma_2 + 9\gamma_3}{4n}\right) \qquad (2.6)$$

## 3. Adapted family of estimators

Following Solanki and Singh (2013), we propose a class of estimators of population mean in simple random sampling using information from a Poisson distributed population, as

$$t_P = \bar{y}_{po} \exp\left(\frac{\alpha(\bar{X} - \bar{x}_{po})}{\bar{X} + \bar{x}_{po}}\right) \qquad (3.1)$$

where $\alpha$ is suitably chosen scalar.

Here, we note that

1. For $\alpha = 0$, $t_P = \bar{y}_{po}$ that is usual unbiased estimator in simple random sampling.
2. For $\alpha = 1$, $t_P = t_{k1}$, that is exponential ratio estimator.
3. For $\alpha = -1$, $t_P = t_{k2}$, that is exponential product estimator.

The proposed class of estimators is a generalized estimator class of unbiased estimators $\bar{y}_{po}$, exponential ratio $t_{k1}$ and exponential product estimator $t_{k2}$

Expressing the estimator $t_1$ in equation (3.1) in terms of e's, we have

$$t_P = \bar{Y}(1 + e_0) \exp\left(\frac{-\alpha e_1}{2 + e_1}\right) \qquad (3.2)$$

Neglecting the terms having power greater than two of the above expression, we have

$$(t_P - \bar{Y}) = \bar{Y}\left[e_0 - \frac{\alpha e_1}{2} - \frac{\alpha e_0 e_1}{2} + \frac{\alpha(\alpha + 2)}{8} e_1^2\right] \qquad (3.3)$$

Taking expectations of both sides of equation (3.3), we get the bias of the estimator $t_1$ to the first degree of approximation, as

$$\text{Bias}(t_P) = \bar{Y}\left[\frac{\alpha(\alpha + 2)}{8} \frac{(\gamma_1 + \gamma_3)}{n} - \frac{\alpha}{2} \frac{\gamma_3}{n}\right] \qquad (3.4)$$

Squaring both sides of equation (3.3) and neglecting terms of e's having power greater than two, we have

$$(t_P - \bar{Y})^2 = \bar{Y}^2\left[e_0^2 + \frac{\alpha^2 e_1^2}{4} - \alpha e_0 e_1\right] \qquad (3.5)$$

Taking expectations of both sides of above expression, we get the MSE of the estimator $t_1$ to the first degree of approximation, as

$$\text{MSE}(t_P) = \frac{\overline{Y}^2}{n}\left[(\gamma_2 + \gamma_3) + \frac{\alpha^2}{4}(\gamma_1 + \gamma_3) - \alpha\gamma_3\right] \quad (3.6)$$

It is interesting to note here that if we put $\alpha = 0, 2, -2, 1$ and $-1$ in equation (3.6), we get respectively the MSE expressions of the estimators $\overline{y}_{po}$, simple ratio, simple product, $t_{k1}$ and $t_{k2}$ up to the first order of approximation.

Differentiating equation (3.6) with respect to $\alpha$ and then equating it to zero, we get the optimum value of $\alpha$ as

$$\alpha^* = \frac{2\gamma_3}{(\gamma_1 + \gamma_3)} \quad (3.7)$$

Putting the optimum value of $\alpha$ into equation (3.6), we get the minimum MSE of the estimator $t_1$ as

$$\text{MSE}(t_P)_{\min} = \frac{\overline{Y}^2}{n}\left[\gamma_2 + \gamma_3\left(1 - \frac{\gamma_3}{(\gamma_1 + \gamma_3)}\right)\right] \quad (3.8)$$

In sampling literature many authors suggested difference type estimators to get more efficient estimates, moving along this direction, we have suggested a difference estimator for $\overline{y}_{po}$ as

$$t_R = \overline{y}_{po} + b(\overline{X} - x_{po}) \quad (3.9)$$

Expressing the estimator $t_R$ in equation (3.9) in terms of e's, we have

$$t_R = \overline{Y}(1 + e_0) - b\overline{X}e_1$$

Or

$$(t_R - \overline{Y})^2 = (\overline{Y}e_0 - b\overline{X}e_1)^2 \quad (3.10)$$

Simplifying and neglecting terms of e's having power greater than two and taking expectations both sides, we have

$$\text{MSE}(t_R) = \frac{1}{n}\left[\overline{Y}^2(\gamma_2 + \gamma_3) + b^2\overline{X}^2(\gamma_1 + \gamma_3) - 2b\overline{Y}\overline{X}\gamma_3\right] \quad (3.11)$$

Differentiating equation (3.11) with respect to $b$ and equating it to zero, we get the optimum value of $b$ as

$$b^* = \frac{\overline{Y}}{\overline{X}} \frac{\gamma_3}{(\gamma_1 + \gamma_3)} \qquad (3.12)$$

Putting the optimum value of b into equation (3.11), we get the minimum MSE of the estimator $t_R$ as

$$MSE(t_R)_{min} = \frac{\overline{Y}^2}{n}\left[\gamma_2 + \gamma_3\left(1 - \frac{\gamma_3}{(\gamma_1 + \gamma_3)}\right)\right] \qquad (3.13)$$

The MSE expression (3.9) is same as the MSE expression of generalised class of exponential estimator defined earlier.

### 4. The Suggested Generalised Class of Estimators

We propose a generalized family of estimators for population mean of the study variable Y, as

$$t_m = \left\{w_1 \overline{y}_{po}\left(\frac{\overline{X}}{\overline{x}_{po}}\right)^\alpha \exp\left(\frac{\eta(\overline{X} - \overline{x}_{po})}{\eta(\overline{X} + \overline{x}_{po}) + 2\theta}\right)\right\} + w_2 \overline{x}_{po} + (1 - w_1 - w_2)\overline{X} \qquad (4.1)$$

where $w_1$ and $w_2$ are suitable constants to be determined such that MSE of $t_m$ is minimum, $\eta$ and $\theta$ are either real numbers or the functions of the known parameters of auxiliary variables such as coefficient of variation $C_x$, skewness $\beta_{1(x)}$, kurtosis $\beta_{2(x)}$ and correlation coefficient $\rho$ (see Sharma and Singh (2013)).

It is to be mentioned that

(i) For $(w_1, w_2) = (1, 0)$, the class of estimator $t_m$ reduces to the class of estimator as

$$t_{mp} = \left\{\overline{y}_{po}\left(\frac{\overline{X}}{\overline{x}_{po}}\right)^\alpha \exp\left(\frac{\eta(\overline{X} - \overline{x}_{po})}{\eta(\overline{X} + \overline{x}_{po}) + 2\theta}\right)\right\} \qquad (4.2)$$

(ii) For $(w_1, w_2) = (w_1, 0)$, the class of estimator $t_m$ reduces to the class of estimator as

$$t_{mq} = \left\{w_1 \overline{y}_{po}\left(\frac{\overline{X}}{\overline{x}_{po}}\right)^\alpha \exp\left(\frac{\eta(\overline{X} - \overline{x}_{po})}{\eta(\overline{X} + \overline{x}_{po}) + 2\theta}\right)\right\} \qquad (4.3)$$

A set of new estimators generated from (4.1) using suitable values of $w_1, w_2, \alpha, \eta$ and $\theta$ are listed in Table 4.1.

**Table 4.1: Set of estimators generated from the class of estimators $t_m$**

| Subset of proposed estimator | $w_1$ | $w_2$ | $\alpha$ | $\eta$ | $\theta$ |
|---|---|---|---|---|---|
| $t_{m1} = \bar{y}_{po}$ | 1 | 0 | 0 | 0 | 1 |
| $t_{m2} = \bar{y}_{po}\left(\dfrac{\bar{X}}{\bar{x}_{po}}\right)$ | 1 | 0 | 1 | 0 | 1 |
| $t_{m3} = \bar{y}_{po}\left(\dfrac{\bar{X}}{\bar{x}_{po}}\right)^{\alpha}$ (Srivastava, 1967) | 1 | 0 | $\alpha$ | 0 | 1 |
| $t_{m4} = \bar{y}_{po}\left(\dfrac{\bar{x}_{po}}{\bar{X}}\right)$ | 1 | 0 | -1 | 0 | 1 |
| $t_{m5} = w_1 \bar{y}_{po}\left(\dfrac{\bar{X}}{\bar{x}_{po}}\right)$ | 1 | 0 | 1 | 0 | 1 |
| $t_{m6} = w_1 \bar{y}_{po}\left(\dfrac{\bar{x}_{po}}{\bar{X}}\right)$ | $w_1$ | 0 | -1 | 0 | 1 |
| $t_{m7} = w_1 \bar{y}_{po}$ | $w_1$ | 0 | 0 | 0 | 1 |

another set of estimators generated from class of estimator $t_{mq}$ given in (4.3) using suitable values of $\eta$ and $\theta$ are summarized in table 4.2

**Table 4.2: Set of estimators generated from the estimator $t_{mq}$**

| Subset of proposed estimator | $\alpha$ | $\eta$ | $\lambda$ |
|---|---|---|---|
| $t_{mq}^{(1)} = \left\{ w_1 \bar{y}_{po}\left(\dfrac{\bar{X}}{\bar{x}_{po}}\right) \exp\left(\dfrac{(\bar{X} - \bar{x}_{po})}{(\bar{X} + \bar{x}_{po}) + 2}\right) \right\}$ | 1 | 1 | 1 |
| $t_{mq}^{(2)} = \left\{ w_1 \bar{y}_{po} \left( \dfrac{\bar{X}}{\bar{x}_{po}} \right) \exp\left( \dfrac{(\bar{X}-\bar{x}_{po})}{(\bar{X}+\bar{x}_{po})+2\rho} \right) \right\}$ | 1 | 1 | $\rho$ |
| $t_{mq}^{(3)} = \left\{ w_1 \bar{y}_{po} \left( \dfrac{\bar{X}}{\bar{x}_{po}} \right) \exp\left( \dfrac{(\bar{X}-\bar{x}_{po})}{(\bar{X}+\bar{x}_{po})+2\bar{X}} \right) \right\}$ | 1 | 1 | $\bar{X}$ |
| $t_{mq}^{(4)} = \left\{ w_1 \bar{y}_{po} \left( \dfrac{\bar{X}}{\bar{x}_{po}} \right) \exp\left( \dfrac{(\bar{X}-\bar{x}_{po})}{(\bar{X}+\bar{x}_{po})} \right) \right\}$ | 1 | 1 | 0 |
| $t_{mq}^{(5)} = \left\{ w_1 \bar{y}_{po} \left( \dfrac{\bar{x}_{po}}{\bar{X}} \right) \exp\left( \dfrac{(\bar{X}-\bar{x}_{po})}{(\bar{X}+\bar{x}_{po})} \right) \right\}$ | -1 | 1 | 1 |
| $t_{mq}^{(6)} = \left\{ w_1 \bar{y}_{po} \left( \dfrac{\bar{X}}{\bar{x}_{po}} \right) \exp\left( \dfrac{\bar{X}(\bar{X}-\bar{x}_{po})}{\bar{X}(\bar{X}+\bar{x}_{po})+2\rho} \right) \right\}$ | 1 | $\bar{X}$ | $\rho$ |
| $t_{mq}^{(7)} = \left\{ w_1 \bar{y}_{po} \exp\left( \dfrac{\bar{X}(\bar{X}-\bar{x}_{po})}{\bar{X}(\bar{X}+\bar{x}_{po})+2\rho} \right) \right\}$ | 0 | $\bar{X}$ | $\rho$ |
| $t_{mq}^{(8)} = \left\{ w_1 \bar{y}_{po} \left( \dfrac{\bar{X}}{\bar{x}_{po}} \right) \exp\left( \dfrac{\rho(\bar{X}-\bar{x}_{po})}{\rho(\bar{X}+\bar{x}_{po})+2\bar{X}} \right) \right\}$ | 1 | $\rho$ | $\bar{X}$ |
| $t_{mq}^{(9)} = \left\{ w_1 \bar{y}_{po} \left( \dfrac{\bar{x}_{po}}{\bar{X}} \right) \exp\left( \dfrac{\rho(\bar{X}-\bar{x}_{po})}{\rho(\bar{X}+\bar{x}_{po})+2\bar{X}} \right) \right\}$ | -1 | $\rho$ | $\bar{X}$ |

---

Expressing (4.1) in terms of e's, we have

$$t_m = w_1 \bar{Y}(1+e_0)(1+e_1)^{-\alpha} \exp\{-ke_1(1+ke_1)^{-1}\} + w_2 \bar{X}(1+e_1) + (1-w_1-w_2)\bar{X}$$

where, $k = \dfrac{\eta}{2(\eta \bar{X} + \theta)}$. (4.4)

Up to the first order of approximation, we have

$$(t_m - \bar{Y}) = [(w_1 - 1)b + w_2 \bar{Y}\{e_0 - ae_1 + de_1^2 - ae_0 e_1\} + w_2 \bar{X} e_1] \qquad (4.5)$$

where $a = (\alpha + k)$, $d = (\bar{Y} - \bar{X})$ and $d = \left\{ \dfrac{3}{2}k^2 + \alpha k + \dfrac{\alpha(\alpha+1)}{2} \right\}$

Squaring both sides of equation (4.5) and neglecting terms of e's having power greater than two, we have

$$(t_m - \overline{Y})^2 = [(1 - 2w_1)d^2 + w_1^2\{d^2 + \overline{Y}^2(e_0^2 + a^2e_1^2 - 2ae_0e_1)\}$$
$$+ w_2^2\overline{X}^2e_1^2 + 2w_1w_2\overline{YX}(e_0e_1 - ae_1^2)] \quad (4.6)$$

Taking expectations both sides, we get the MSE of the estimator $t_m$ to the first order of approximation as

$$MSE(t_m) = [(1 - 2w_1)d^2 + w_1^2 A + w_2^2 B + 2w_1w_2 C] \quad (4.7)$$

where,

$$A = d^2 + \frac{\overline{Y}^2}{n}(\lambda_2 + a^2\lambda_1 - 2a\gamma_2),$$

$$B = \overline{X}^2 \frac{\lambda_1}{n},$$

$$C = \overline{YX}(\gamma_2 - a\lambda_1).$$

The optimum values of $w_1$ and $w_2$ are obtained by minimizing (4.7) and is given by

$$w_1^* = \frac{d^2 B}{(AB - C^2)} \quad \text{And} \quad w_2^* = \frac{-d^2 C}{(AB - C^2)} \quad (4.8)$$

Substituting the optimal values of $w_1$ and $w_2$ in equation (4.7) we obtain the minimum MSE of the estimator $t_m$ as

$$MSE_{min}(t_m) = d^2 \left[1 - \frac{d^2 B}{(AB - C^2)}\right] \quad (4.9)$$

Putting the values of A, B, C, d and simplifying, we get the minimum MSE of estimator $t_m$ as

$$MSE_{min}(t_m) = \left[\frac{\lambda_2^2(\lambda_1\lambda_2 - \gamma_2^2)}{\left[\lambda_1 + \frac{\lambda_2^2}{d^2 n}(\lambda_1\lambda_2 - \gamma_2^2)\right]}\right] \quad (4.10)$$

## 5. Efficiency Comparisons

First we compare the efficiency of the proposed estimator $t_p$ with Koyuncu and Ozel (2013) exponential ratio estimator,

$$\text{MSE}(t_{k1}) - \text{MSE}(t_P)_{\min} \geq 0$$

If

$$\overline{Y}^2\left(\frac{\gamma_1 + 4\gamma_2 + \gamma_3}{4n}\right) - \frac{\overline{Y}^2}{n}\left[\gamma_2 + \gamma_3\left(1 - \frac{\gamma_3}{(\gamma_1 + \gamma_3)}\right)\right] \geq 0$$

or

$$(\gamma_1 + \gamma_3 + 4\gamma_2) \geq 4\gamma_2 + 4\gamma_3\left(1 - \frac{\gamma_3}{(\gamma_1 + \gamma_3)}\right)$$

or

$$(\gamma_1 + \gamma_3)^2 \geq 4\gamma_1\gamma_3 \tag{5.1}$$

We observe that the condition listed in (5.1) shows that proposed class of estimators is always better than the estimator of Koyuncu and Ozel (2013).

Next, we compare the adapted Solanki and Singh (2013) estimator with proposed family of estimators $t_p$,

$$\text{MSE}(t_{k1}) - \text{MSE}(t_P)_{\min} > 0$$

If

$$\overline{Y}^2\left(\frac{\gamma_1 + 4\gamma_2 + 9\gamma_3}{4n}\right) \geq \frac{\overline{Y}^2}{n}\left[\gamma_2 + \gamma_3\left(1 - \frac{\gamma_3}{(\gamma_1 + \gamma_3)}\right)\right] \tag{5.2}$$

Or

$$(\gamma_1 + 5\gamma_3) \geq -\frac{4\gamma_3^2}{(\gamma_1 + \gamma_3)} \tag{5.3}$$

On solving we observe that above conditions holds always true.

Finally, we are comparing the mean square errors of generalised class of estimators and exponential class of estimators

$$\text{MSE}(t_P)_{min} - \text{MSE}(t_m)_{min} \geq 0$$

If

$$\frac{\overline{Y}^2}{n}\left[\gamma_2 + \gamma_3\left(1 - \frac{\gamma_3}{(\gamma_1 + \gamma_3)}\right)\right] - \left[\frac{\lambda_2^2(\lambda_1\lambda_2 - \gamma_2^2)}{\left[\lambda_1 + \frac{\lambda_2^2}{d^2n}(\lambda_1\lambda_2 - \gamma_2^2)\right]}\right] \geq 0 \quad (5.4)$$

or

$$\lambda_2^2(\lambda_1\lambda_2 - \gamma_2^2) \geq (n-1)\lambda_1 d^2 n \quad (5.5)$$

for the condition given in (5.5) the generalised class of estimators performs well.

## 6. Empirical study

**Data Statistics:** To illustrate the efficiency of proposed estimators in the application, we consider the earthquake data of Turkey for the numerical comparisons of the proposed estimators and existing exponential ratio and product estimators in the simple random sampling. The data is obtained from the data base of Kandilli Observatory, Turkey. Earthquake is an unavoidable natural disaster for Turkey since a significant portion of turkey is subject to frequent destructive mainshocks, their foreshocks and aftershocks sequence. We consider the mainshocks that occurred in 1900 and 2011 having surface wave magnitudes $M_S \geq 0.5$, their foreshocks within five days week with $M_S \geq 0.3$ and aftershocks within one month with $M_S \geq 0.4$. In this area, 109 mainshocks with surface magnitude $M_S \geq 0.5$ have occurred between 1900 and 2011. The population consists of the destructive earthquakes. In the population data set the number of aftershocks is a study variable and the number of foreshocks is an auxiliary variable. Note that we take sample size n=20. The MSE values of the proposed estimators are computed with considering the distribution of study and auxiliary variables. To obtain the distribution of these variables, we fit the Poisson distribution to the earthquake data set. To obtain the $\rho_{x_{po}y_{po}}$ for the Poisson distributed data,

Turkey is divided into three main neotectonic domains based on the neotectonic zones of Turkey. The foreshocks in Turkey are separated according to these neotectonic zones. In this way, the parameters $\gamma_1$, $\gamma_2$ and $\gamma_3$ are obtained. According to the goodness of the fit test, it is obvious that the Poisson distribution fits the number of shocks for Region 1 with parameter $\gamma_1 = 4.1813$ $\left(\chi^2 = 0.048, p-\text{value} = 0.043\right)$ and $\gamma_2 = 8.104$ $\left(\chi^2 = 0.014, p-\text{value} = 0.032\right)$ for Region 2, and $\gamma_3 = 2.112$ $\left(\chi^2 = 0.013, p-\text{value} = 0.025\right)$ for Region 3. Then the correlation between the study variable and auxiliary variable is positive $\rho_{x_{po}y_{po}} = 0.712$ and it can be noted that the number of foreshocks is related to the number of aftershocks therefore, ratio estimators are preferable in this case.

We compute the percent relative efficiency (PRE) of different estimators using given expression

$$\text{PRE}(\bullet) = \frac{\text{Var}(\bar{y}_{po})}{\text{MSE}(\bullet)} \times 100$$

**Table 6.1: PRE values of Estimators**

| Estimators | PREs |
|---|---|
| $\bar{y}_{po}$ | 100.00 |
| $t_r$ | 100.00 |
| $t_{k1}$ | 103.25 |
| $t_{k2}$ | 72.31 |
| $t_P$ | 106.73 |
| $t_R$ | 106.73 |
| $t_m$ | **9937.42** |

Table 6.2 exhibits that the percent relative efficiency of the existing and proposed generalised class of estimators using the Poisson distributed population. From the above table we analyse that the usual unbiased estimator and ratio estimator $t_r$ has same PRE's whereas exponential ratio and exponential product estimators are more and less efficient than the simple ratio estimator. The exponential product estimator $t_{k2}$ is less efficient than all the other estimators since the correlation is positive and high. It can also be seen that the proposed generalized exponential estimator $t_P$ and difference estimator $t_R$ are equally efficient but more efficient than the ratio estimator $t_r$, exponential ratio and exponential product estimators (due to Koyuncu and Ozel, 2013). The generalised class of estimators $t_m$ has maximum PRE among all the estimators considered here therefore it is best estimator in the sense of having the least MSE. Most of the members of the generalised class of estimators $t_m$ have approximately same or less efficient than the minimum mean square error of the estimator $t_m$ of the proposed class of estimators therefore it is worthless to mention their mean square errors here though the members are shown in the table 4.1 and 4.2.

**Conclusion**

In this study we have proposed an exponential class of estimators, a difference estimator and a generalised class of estimators considering the distribution of the study and auxiliary variable. The Bias and mean square errors expressions are derived up to the first order of approximation. It has been shown under theoretical and empirical comparisons that the two proposed estimators exponential class of estimators and difference estimator are equally efficient but always more efficient than all other available estimators in this paper. It is also observed that the generalised class of estimators given in section 4 are always more efficient than the other estimators using Poisson distributed population. The proposed generalized class of estimators is advantageous in the sense that the properties of the estimators, which are members of the proposed class of estimators, can be easily obtained from the properties of the proposed generalized class. Thus the study unifies properties of several estimators for population mean using Poisson distributed population